\documentclass[11pt,a4paper,srcltx]{article}
\usepackage{amsmath}
\usepackage[colorlinks=true]{hyperref}
\usepackage{xcolor}
\usepackage[overload]{empheq}

\def\Dom{{\mathrm{Dom}\,}}

\def\cD{{\cal D}}

\def\cE{{\cal E}}

\def\cT{{\cal T}}

\def\rank{{\hbox{\rm Rank}}}
\def\cX{{\cal X}}

\def\Argmin{{\mathrm{Argmin}}}

\newcommand{\half}{ \mbox{\small$\frac{1}{2}$}}

\usepackage{amsfonts}
\usepackage{epsfig}

\def\Ker{\mathop{\hbox{\rm Ker}}}

\usepackage{amssymb}

\usepackage{graphicx}
\def\bR{{\mathbf{R}}}

\def\bS{{\mathbf{S}}}

\newtheorem{lemma}{Lemma}[section]
\newtheorem{corollary}{Corollary}[section]
\newtheorem{proposition}{Proposition}[section]

\newtheorem{theorem}{Theorem}[section]

\def\qed{$\Box$}

\def\mfn{{\mathfrak{n}}}
\def\abs{{\mathrm{abs}}}

\newcommand{\aic}[2]{{  #2}}

\newcommand{\be}{\begin{eqnarray}}
\newcommand{\ee}[1]{\label{#1}\end{eqnarray}}
\newcommand{\nn}{\nonumber \\}
\newcommand{\ese}{\end{eqnarray*}}
\newcommand{\bse}{\begin{eqnarray*}}
\newcommand{\rf}[1]{~(\ref{#1})}

\newcommand{\hide}[1]{{}}
\def\qed{\hfill$\Box$\par}

\def\qed{\hfill$\Box$\par}

\def\mfn{{\mathfrak{n}}}
\title{{Aggregating regular norms}}

\author{
Anatoli Juditsky
\footnotemark[1]%\ \footnotemark[5]
\and Arkadi Nemirovski\footnotemark[2]}
\begin{document}\maketitle
\renewcommand{\thefootnote}{\fnsymbol{footnote}}
\footnotetext[1]{LJK, Universit\'e Grenoble Alpes, 700 Avenue Centrale,  38401 Domaine Universitaire de Saint-Martin-d'Hères, France,
{\tt anatoli.juditsky@univ-grenoble-alpes.fr}}
\footnotetext[2]{Georgia Institute of Technology, Atlanta, Georgia 30332, USA, {\tt nemirovs@isye.gatech.edu}\\
\aic{Research of this author}{This work} was supported by MIAI {@} Grenoble Alpes (ANR-19-P3IA-0003).}
\renewcommand{\thefootnote}{\arabic{footnote}}
\date{}

\begin{abstract}
The subject of this paper is regularity-preserving aggregation of {\sl regular} norms on finite-dimensi\-o\-nal linear spaces. Regular norms were introduced
in \cite{JudNem08} and are closely related to ``type 2'' \aic{norms}{spaces}  \cite[Chapter 9]{Milman} playing important role in 1) high-dimensional convex geometry and probability in Banach spaces, see \cite{Pis75,Pin92,Pin94,Milman,Talagrand}, and in 2) design of proximal first order algorithms for large-scale convex optimization with dimension-independent, or nearly so, complexity. Regularity, with moderate parameters, of a norm makes applicable, in a dimension-independent fashion, \aic{numerous}{various}  geometric, probabilistic, and optimization-related results, thus motivating the subject of this paper---aggregations of regular norms resulting in  controlled (and moderate) {inflation} of regularity parameters.
\end{abstract}

\section{Introduction}
This paper focuses on {\sl regular norms} on (finite-dimensional) linear spaces---the notion introduced in \cite{JudNem08} in connection with developing dimension-independent bounds on large deviations of sums of i.i.d. random vectors or vector-valued matringales in normed spaces.  Postponing definition of regularity till Section \ref{preliminaries}, right now it suffices to say that this property is quantified by a pair of reals $(\varkappa\geq1,\varsigma\geq1)$, the \aic{less}{smaller} the better,  and is one of numerous existing ``measures of deviation'' of the norm in question from the standard Euclidean norm (for the latter, $\varkappa=\varsigma=1$). \par
There are two broad application areas where regularity of norms is of significant interest. The first is high-dimensional convex geometry and probability in Banach spaces, see \cite{Pis75,Pin92,Pin94,Milman,Talagrand} and references therein. Numerous results in these applications
assume that the normed space in question is of ``type $p$'' (see \cite[Chapter 9]{Milman}), the case of $p=2$ being, in a sense, the best. The qualitative property of a normed space $X$ to be of type 2 (which is trivially true for finite-dimensional spaces we are working with) has numerical \aic{characterization,}{characteristics} $T(X)$ {which enters} the quantitative components of the associated geometric and probabilistic results, the smaller is this constant, the better. It is immediate to see that when the norm on $X$ is $(\varkappa,\varsigma)$-regular,one has $T(X)\leq\varsigma^2\sqrt{\varkappa}$. That is, regularity of the norm in question with ``moderate'' constant makes applicable, {\sl in a dimension-independent fashion}, numerous powerful geometric and probabilistic results. Another application domain is that of first order proximal type algorithms for large-scale deterministic and stochastic
convex optimization (see, e.g., \cite{Beck,LMCO,Lan,Yura}). It turns out that regularity of a norm allows to equip the ball of the {\sl conjugate} norm with a ``good proximal setup'' resulting in {\sl dimension-independent} iteration complexity of associated proximal algorithms for minimizing convex functions over these balls.
\par
Basic examples of regular norms (\cite{JudNem08}) are:
\begin{enumerate}
\item $\|\cdot\|_2$ on $\bR^n$ is $(1,1)$ regular. The norm $\|\cdot\|_p$ on $\bR^n$, $2\leq p\leq \infty$, is $(\varkappa,\varsigma)$-regular with $\varkappa=O(1)\min[p,\ln(n+1)]$ and $\varsigma=O(1)$ (here and in what follows $O(1)$ stands for an absolute constant).
\item When $2\leq p\leq\infty$, the Schatten norm $\|\sigma(x)\|_p$ on the space $\bR^{m\times n}$ of $m\times n$ matrices, $\sigma(x)$ being the singular spectrum of matrix $x$, is $(\varkappa,\varsigma)$-regular with $\varkappa=O(1)\min[p,\ln(\min(m,n))+1]$ and $\varsigma=O(1)$.
\end{enumerate}
It should be stressed that parameters of regularity in these examples are either dimension-independent, or deteriorate only logarithmically as the dimension grows.
\par
The potential applications of regular norms outlined above motivate the goal of this paper---{\sl investigating regularity-preserving
``aggregation'' of regular norms.}
\par
The main body of the paper is organized as follows:
{We recall the notion of regular norm in Section \ref{regnorms}  and then present our principal contributions---Theorems \ref{the1Gen} and \ref{theass}---in Section \ref{sec:regaggn}; preliminary versions of these results can be found in our preprint \cite{JudNem08}.  Then in Section \ref{sec:elln} we apply these results to prove the regularity of {\em ellitopic and spectratopic norms,} as defined in \cite{JudNem20,JKN21}. All proofs are postponed to the appendix.}

\section{Problem statement and main results}\label{preliminaries}
\subsection{Regular norms}\label{regnorms} Let $\varkappa\geq1$ and $\varsigma\geq 1$ be two reals. We refer to the norm $\mfn(\cdot)$ on $\bR^n$ as {\em $\varkappa$-smooth} if the function
$\Phi(x)=\mfn^2(x)$ is continuously differentiable and satisfies the relation
$$
\forall x,h: \Phi(x+h)\leq \Phi(x)+[\Phi'(x)]^Th+\varkappa \Phi(h),
$$
or, which is the same (see Lemma \ref{simplem}), the gradient of $\Phi$ is Lipschitz continuous with Lipschitz constant $2\varkappa$:
$$
\mfn_*(\nabla \Phi(x)-\nabla \Phi(y))\leq 2\varkappa \mfn(x-y)
$$
where $\mfn_*(\cdot)$ is the norm conjugate to $\mfn(\cdot)$.
\par
Following \cite{JudNem08}, we say that a norm $\|\cdot\|$ on $\bR^n$ is {\em $(\varkappa,\varsigma)$-regular} if there exists a $\varkappa$-smooth norm $\mfn(\cdot)$ on $\bR^n$ such that $\|\cdot\|$ is within factor $\varsigma$ from $\mfn(\cdot)$:
$$
\varsigma^{-1}\mfn(\cdot)\leq \|\cdot\|\leq \varsigma \mfn(\cdot).
$$
\par Clearly, if $\|\cdot\|$ is $\varkappa$-smooth (or $(\varkappa,\varsigma)$-regular) and $x\mapsto
 Ax$ is a linear embedding of $\bR^m$ into $\bR^n$ (so that $\Ker A=\{0\}$), then the norm
$\|y\|_A=\|Ay\|$ on $\bR^m$ is  $\varkappa$-smooth  (resp., $(\varkappa,\varsigma)$-regular). Besides this, if $\|\cdot\|'$ is within factor $\alpha\geq1$ from $(\varkappa,\varsigma)$-regular norm $\|\cdot\|$, that is,
$$
\alpha^{-1}\|\cdot\|'\leq \|\cdot\|\leq\alpha\|\cdot\|',
$$
then $\|\cdot\|'$ is $(\varkappa,\alpha\varsigma)$-regular.
\par
The examples of regular norms presented in the introduction stem straightforwardly from the following facts (\cite{JudNem08}; in proximal minimization context, these facts are reproduced in \cite{WithYu,WithRo})  which \aic{in what follows are important by their own right}{are important in what follows and in their own right}:
{\sl when $2\leq p<\infty$,}
\begin{itemize}
\item%[{\bf A}]
{\sl The norm $\|\cdot\|_p$ on $\bR^n$ is $(p-1)$-smooth }\cite[Example 3.2, Section 4.1.1]{JudNem08}.
\item%[{\bf B}]
{\sl The Schatten norm $\|\sigma(X)\|_p$ on the space $\bR^{m\times n}$ of $m\times n$ matrices ($\sigma(x)\in\bR^{\min[m,n]}$  stands for the  singular spectrum of $X$) is $\max[2,p-1]$-smooth} \cite[Example 3.3, Section 4.1.1]{JudNem08}.
\end{itemize}

\subsection{Main results}
\label{sec:regaggn}
Our main results on aggregation of regular norms are as follows:
\begin{theorem}\label{the1Gen}
Let $\theta(\cdot):\bR^K_+\to\bR$ be a convex continuous homogeneous, of degree 1, function which is monotone on $\bR^K_+$  and positive outside of  the origin. Let $\|\cdot\|_i$ be $(\varkappa,\varsigma)$-regular norms on $\bR^{n_i}$, $1\leq i\leq K$.
 Then the aggregated norm
 \begin{equation}\label{aggregatednorm}
 \|[x_1;...;x_K]\|=\theta^{1/2}(\|x_1\|_1^2,...,\|x_K\|_K^2)
\end{equation}
(clearly, this indeed is a norm on $\bR^{n_1+...+n_K}$) is $(c_1[\ln(K+1)+\varkappa],c_2\varsigma)$-regular, for properly selected absolute constants $c_1,c_2$.
\end{theorem}

Our next result is the refinement of Theorem \ref{the1Gen} in the situation where the ``aggregating function'' $\theta$ is a regular absolute norm:
\begin{theorem}\label{theass} Let $\theta$ be an absolute norm on $\bR^K$ which is $(\overline{\varkappa},\overline{\varsigma})$-regular, and let $\|\cdot\|_i$ be $(\varkappa',\varsigma')$-regular norms on $\bR^{n_i}$, $i\leq K$. Then the norm
\[
 \|[x_1;...;x_K]\|=\theta^{1/2}(\|x_1\|_1^2,...,\|x_K\|_K^2)
 \]
is $(\varkappa,\varsigma)$-regular on $\bR^{n_1+...+n_K}$, with
$
\varkappa =2\overline{\varkappa}+\varkappa'$ and $\varsigma=\varsigma'\sqrt{\overline{\varsigma}}.$
\end{theorem}
Finally, we have
\begin{proposition}\label{thefact} Let $\|\cdot\|$ be a $(\varkappa,\varsigma)$-regular norm on $\bR^n$ and $x\mapsto Px: \bR^n\to\bR^p$ be an onto mapping. Then the factor-norm
$$
\|u\|'=\min_x\{\|x\|:Px=u\}
$$
induced on $\bR^p$ by $\|\cdot\|$, $P$, is $(\varkappa,\varsigma)$-regular.
\end{proposition}
\paragraph{\bf Remark:} Theorems \ref{the1Gen}, \ref{theass} and Proposition \ref{thefact} \aic{look}{are formulated} as existence results stating that under such and such assumptions such and such norm $\|\cdot\|$ is regular with such and such parameters $\varkappa$, $\varsigma$, so that there exists a $\varkappa$-smooth norm $\|\cdot\|'$ which is within factor $\varsigma$ from $\|\cdot\|$. However, in some applications, e.g., in proximal-type minimization, ``simple {\sl existence}'' of $\|\cdot\|'$ is not enough---we should know what this norm is. In this respect, it should be stressed that, as is seen from the proofs, the above statements are ``constructive.'' For example, in the context of Theorem \ref{the1Gen}, given $\kappa$-smooth norms which are within the factor $\varsigma$ from the norms $\|\cdot\|_i$, we know how to convert these norms and $\theta(\cdot)$ into a $(c_1[\ln(K+1)+\varkappa])$-smooth norm which is within factor $c_2\varsigma$ of the aggregated norm (\ref{aggregatednorm}), and similarly for Theorem \ref{theass} and Proposition \ref{thefact}.

\section{Illustration: regularity of ellitopic and spectratopic norms}
\label{sec:elln}
\subsection{Ellitopes and spectratopes}
\paragraph{Ellitopes.} A {\em basic ellitope} \cite[Section 4.2.1]{JudNem20} in $\bR^n$ is a set $\cX$ given as
\begin{equation}\label{eeq1}
\cX=\{x\in\bR^n:\exists t\in\cT:x^TT_ix\leq t_i,\,1\leq i\leq K\},
\end{equation}
where $T_i\succeq0$, $\sum_iT_i\succ0$, and $\cT$ is a convex compact set in $\bR^K_+$ which contains a positive vector and is monotone: whenever $0\leq t'\leq t\in\cT$, we have $t'\in\cT$. Here $A\succeq B\Leftrightarrow B\preceq A$ ($A\succ B\Leftrightarrow B\prec A$)   mean that $A,B$ are symmetric matrices of the same size such that $B-A$ is positive semidefinite (resp., positive definite).
\par
{\sl Ellitope} $\overline{X}\subset\bR^p$ is a linear image of basic ellitope:
\begin{equation}\label{eeq100}
\overline{\cX}=P\cX\hbox{\ with $\cX$ given by (\ref{eeq1})},
\end{equation}
where $P\in\bR^{p\times n}$. \par
The simplest examples of basic ellitopes are
\begin{itemize}
\item Bounded intersections $\cX=\{x:x^TT_kx\leq1,k\leq K\}$, $T_k\succeq0$,  of centered at the origin ellipsoids/elliptic cylinders.
\item $\|\cdot\|_p$-balls, $p\geq 2$: $\{x\in\bR^n:\|x\|_p\leq1\}=\{x\in\bR^n:\exists t\in\cT:x_i^2\leq t_i,\,i\leq n\}$, $\cT=\{t\in\bR^n_+:\|t\|_{p/2}\leq1\}$.
\end{itemize}
\paragraph{Spectratopes.} A {\sl basic spectratope} \cite[Section 4.3.1]{JudNem20} in $\bR^n$ is a set $\cX$ given as
\begin{equation}\label{eeq143}
\cX=\{x\in\bR^n:\exists t\in\cT:S_i^2[x]\preceq t_iI_{d_i},\,1\leq i\leq K\},
\end{equation}
where $S_i[x]=\sum_{j=1}^nx_jS^{ij}$, $S^{ij}\in\bS^{d_i}$, are linear mappings taking values in the spaces $\bS^{d_i}$ of $d_i\times d_i$ symmetric matrices and such that $S_i[x]=0$ for all $i\leq K$ if and only if $x=0$, and $\cT\subset\bR^K_+$ is of the same type as in the definition of a basic ellitope.
\par
{\sl Spectratope} $\overline{X}\subset\bR^p$ is a linear image of basic spectratope:
\begin{equation}\label{eeq10043}
\overline{\cX}=P\cX\hbox{\ with $\cX$ given by (\ref{eeq143})},
\end{equation}
where $P\in\bR^{p\times n}$.
\par
The simplest examples of basic spectratopes are
\begin{itemize}
\item The unit ball of spectral norm $\|\cdot\|_{2,2}$ in the space $\bS^m$ of symmetric $m\times m$ matrices:
$$
\{x\in\bS^m:\|x\|_{2,2}\leq1\}=\{x\in \bS^m: \exists t\in\cT:=[0,1]: S^2[x]\preceq tI_m\},\,\,S[x]\equiv x.
$$
\item The unit ball of the spectral norm $\|\cdot\|_{2,2}$ in the space $\bR^{m\times n}$ of $m\times n$ matrices:
$$
\{x\in\bR^{m\times n}:\|x\|_{2,2}\leq1\}=\{x\in \bR^{m\times n}:\exists t\in\cT:=[0,1]: S^2[x]\preceq tI_{m+n}\},\,S[x]=\left[\begin{array}{c|c}&x\cr\hline x^T&\cr\end{array}\right].
$$
\end{itemize}
\paragraph{Calculus of ellitopes/spectratopes.}
Families of ellitopes and of spectratopes are rather rich: as shown in \cite[Section 4.6]{JudNem20}, each family is closed w.r.t. taking finite intersections, direct products, arithmetic sums, images under linear mappings and inverse images under linear embeddings. Furthremore, ellitopic description \rf{eeq1},\rf{eeq100} and spectratopic description \rf{eeq143}, \rf{eeq10043} of the results of these operations are readily given by ellitopic and, respectively, spectratopic descriptions of the operands. It should also be mentioned that every ellitope is a spectratope as well.
\paragraph{Ellitopic and spectratopic norms.} The sets $\cT\subset\bR^K_+$  participating in the definition of ellitopes and spectratopes are exactly sets of the form
\begin{equation}\label{theta}
\cT=\{t\in\bR^K_+: \theta(t)\leq1\}
\end{equation}
stemming from functions $\theta(\cdot):\bR^K_+\to\bR$ which are
\begin{enumerate}
\item convex, continuous, positive outside of  the origin, and positively homogeneous, of homogeneity degree 1,
\item monotone on $\bR^K_+$: $0\leq t\leq t'\Rightarrow \theta(t)\leq\theta(t').$
\end{enumerate}
Next, basic ellitopes and spectratopes in $\bR^n$ are convex compact sets which are symmetric w.r.t. the origin and have the origin in their interior. Consequently, such a set $\cX$ is the unit ball of a norm $\|\cdot\|_\cX$ on $\bR^N$; this norm is

%\begin{equation}\label{thenormis}
%\|x\|_\cX=\left\{\begin{array}{ll}\theta^{1/2}\Big([\|T_1^{1/2}x\|_2^2;...;\|T_K^{1/2}x\|_2^2]\Big),&\hbox{when $\cX$ is the ellitope (\ref{eeq1})}\\
%\theta^{1/2}\Big([\|S_1[x]\|_{2,2}^2;...;\|S_K[x]\|_{2,2}^2]\Big),&\hbox{when $\cX$ is the spectratope (\ref{eeq143})}\\
%\end{array}\right.,
%\end{equation}
\begin{subequations}
\label{eq:thenormis-all}
\begin{empheq}[left={\|x\|_\cX=\empheqlbrace\,}]{align}
&\theta^{1/2}\Big([\|T_1^{1/2}x\|_2^2;...;\|T_K^{1/2}x\|_2^2]\Big),\quad\hbox{when $\cX$ is the ellitope (\ref{eeq1})}\label{thenormis.a}\\
&\theta^{1/2}\Big([\|S_1[x]\|_{2,2}^2;...;\|S_K[x]\|_{2,2}^2]\Big),\quad\hbox{when $\cX$ is the spectratope (\ref{eeq143})}\label{thenormis.b}
\end{empheq}
\end{subequations}
where $\theta(\cdot)$ is given by (\ref{theta}).
\par
Finally, if $\overline{\cX}$ is an ellitope (a spectratope) given by (\ref{eeq100}) (resp., by (\ref{eeq10043})) {\sl and the corresponding mapping $x\mapsto Px$ is an onto one}, $\overline{\cX}$ also is the unit ball of certain norm
$\|\cdot\|_{\overline{\cX}}$ which is nothing but the factor-norm induced by $\cX$ and $P$:
\begin{equation}\label{factornorm}
\|u\|_{\overline{\cX}}=\min_x\left\{\|x\|_\cX: Px=u\right\}.
\end{equation}

\paragraph{The result} of this section is as follows:
\begin{theorem}\label{theellispec}
There exist absolute constants $c_1$, $c_2$ such that
\item[(i)] Whenever $\overline{\cX}$ is given by {\rm  (\ref{eeq1}), (\ref{eeq100})}, and {\rm (\ref{theta})} with $\rank(P)=p$,  the ellitopic norm $\|\cdot\|_{\overline{\cX}}$ given by {\rm \rf{thenormis.a}} and {\rm (\ref{factornorm})} is $(c_1\ln(K+1),c_2)$-regular.
\item[(ii)] Whenever $\overline{\cX}$ is given by {\rm (\ref{eeq143}), (\ref{eeq10043}), and (\ref{theta})} with $\rank(P)=p$,  the spectratopic norm $\|\cdot\|_{\overline{\cX}}$ given by {\rm \rf{thenormis.b}} and {\rm (\ref{factornorm})} is $(c_1[\ln(K+1)+\max_{i\leq K}\ln(d_i)],c_2)$-regular.
\end{theorem}
%{\bf Proof.} {\em (i):} By Theorem \ref{the1Gen} as applied to the (1,1)-regular norms $\|\cdot\|_i=\|\cdot\|_2$ on $\bR^n$, the norm $\|[x_1;...;x_K]\|:=\theta^{1/2}(\|x_1\|_2^2,...,\|x_K\|_2^2)$ on $\bR^{nK}$ is $(O(1)\ln(K+1),O(1))$-regular. With $\cX$ given by (\ref{eeq1}), the norm $\|x\|_+:=\theta^{1/2}(\|T_1^{1/2}x\|_2^2,...,\|T_k^{1/2}x\|_2^2)$ is the superposition of $\|\cdot\|$ and linear embedding $x\mapsto[T_1^{1/2}x;... ;T_k^{1/2}x]$ and therefore is regular with the same parameters as $\|\cdot\|$, see Section \ref{regnorms}. Finally,  $\|\cdot\|_{\overline{\cX}}$ is the factor-norm of $\|\cdot\|_+$, and passing from a norm to its factor-norm preserves regularity parameters by Proposition \ref{thefact}.
%\par
%{\em (ii):} The norms $\|\cdot\|_i=\|\cdot\|_{2,2}$ on $\bS^{d_i}$ are $(O(1)\ln(d_i+1),O(1))$-regular, see ``basic facts'' at the end of in Section \ref{regnorms}.
%By Theorem \ref{the1Gen}) as applied to these norms, the norm $\|[y_1;...;y_K]\|:=\theta^{1/2}(\|y_1\|_{2,2}^2,...,\|y_K\|_{2,2}^2)$ on $\bS^{d_1}\times...\times\bS^{d_K}$ is $(O(1)[\ln(K+1)+\max_{i\leq K}\ln(d_i)],O(1))$-regular.
%Similarly to the case of $(i)$, with $\cX$ given by \aic{(\ref{eeq100})}{\rf{eeq10043}}, the norm $\|\cdot\|_\aic{\cX}{+}$ \aic{}{which}
%is the superposition of $\|\cdot\|$ and linear embedding $x\mapsto(S_1(x)_,...,S_K(x))$ \aic{and  therefore}{} is regular with the same parameters as $\|\cdot\|$, and these parameters are inherited by the factor-norm $\|\cdot\|_{\overline{\cX}}$. \qed
\appendix
\section{Appendix}
\subsection{Proof of Theorem \ref{the1Gen}}
\subsubsection{Preliminary step}\label{prelim}
\aic{We start with demonstrating}{Let us show first} that in order to prove Theorem \ref{the1Gen} ``as is,'' it suffices to prove it in the case where function $\theta$, in addition to
the properties postulated in theorem's premise (which are convexity and continuity on $\bR^K_+$, homogeneity of degree 1, monotonicity, and positivity outside of the origin), possesses two additional properties, namely,
\begin{enumerate}
\item[(a)] is continuously differentiable outside of the origin,
and
\item[(b)] satisfies the relation
\begin{equation}\label{eeq334}
t>0\Rightarrow 1/K\leq[\nabla\theta(t)]_i\leq (1+1/K),\,1\leq i\leq K.
\end{equation}
\end{enumerate}
Indeed,
\paragraph{1$^o$.} Let $\mfn(t)=\theta(\abs[t])$, where $\abs[[t_1;...;t_K]]=[|t_1|;...;|t_K|]$, so that $\mfn(t)$ is an absolute norm on $\bR^K$. Denoting by $\mfn_*(\cdot)$ the Fenchel conjugate of $\mfn(\cdot)$, the Fenchel conjugate $\vartheta_\epsilon$ of the function $\tfrac{1}{2}\mfn_*^2(\cdot)+\tfrac{\epsilon}{2}\|\cdot\|_2^2$
is a function of the form $\tfrac{1}{2}\theta_\epsilon^2(\abs[\cdot])$  where
$\theta_\epsilon:\,\bR^K_+\to\bR$ is convex homogeneous of degree 1 monotone and positive outside of the origin with Lipschitz continuous gradient; along with $\vartheta_\epsilon$, $\theta_\epsilon$ is continuously differentiable outside of the origin. Taking $\epsilon>0$ small enough, we can make $\theta_\epsilon(\cdot)$ to be  within the factor, say, 1.1, from $\theta(\cdot)$. It clearlty suffices to prove Theorem \ref{the1Gen} with $\theta_\epsilon$ in the role of $\theta$, so that from now on let us assume that $\theta$, in addition to what is postulated in the premise of Theorem \ref{the1Gen}, satisfies (a).
\paragraph{2$^o$.} Now let $s_i=\max \{t_i:\,t\geq0,\,\theta(t)\leq1\}$. By diagonal scaling of $\bR^K$---passing from $\theta(t)$ to $\theta(\cD t)$ with positive definite diagonal $\cD$---we can assume that $s_i=1$, and such scaling does not affect the statement we want to prove. Thus, in addition to what has been already assumed, from now on we assume that $s_i:=\max_t\{t_i:\theta(t)\leq1\}=1$ for all $i$. Monotonicity of $\theta(\cdot)$ implies that the basis orths $e_i$ belong to the set $\Theta=\{t\geq0:\theta(t)\leq1\}$, whence $\theta(t)\leq\sum_it_i$, $t\geq0$. Besides this, under the assumption that $s_i=1$, $i\leq K$, $\Theta$ is contained in the box  $\{t:\,0\leq t_i\leq1,\,i\leq K\}$, and thus is contained in the set $\{t\geq0:\sum_it_i\leq K\}$, implying that $\theta(t)\geq K^{-1}\sum_it_i$, $t\geq0$.
Now let $\overline{\theta}(t)=\theta(t)+K^{-1}\sum_it_i$; by the above, $\theta(\cdot)\leq \overline{\theta}(\cdot)\leq 2\theta(\cdot)$. At the same time, for $t\geq0$, $t\neq0$, one has $[\nabla\overline{\theta}(t)]_i\geq K^{-1}$ for all $i$. Aside of this, we have $\overline{\theta}(t)\leq (1+1/K)\sum_it_i$ for all $t\geq0$, implying that for all $s\geq0$
$$
(1+1/K)[\sum_jt_j+s]\geq \overline{\theta}(t+se_i)\geq \overline{\theta}(t)+s[\nabla\overline{\theta}(t)]_i.
$$
We conclude that $[\nabla\overline{\theta}(t)]_i\leq 1+1/K$ for all $i$ and all nonzero $t\geq0$. The bottom line is that $\overline{\theta}$ is shares all properties mentioned in the beginning of Section \ref{prelim}, while $\theta(\cdot)$ is within factor 2 from $\overline{\theta}(\cdot)$. As a result, the norms ${\theta(\|x_1\|_1^2,...,\|x_K\|_K^2)^{1/2}}$ and  ${\overline{\theta}(\|x_1\|_1^2,...,\|x_K\|_K^2)^{1/2})^{1/2}}$ are within factor $\sqrt{2}$ from each other. We conclude that in order to prove Theorem \ref{the1Gen} it suffices to prove it in the case when $\theta$, in addition to what is stated in the premise of the theorem, possesses properties (a) and (b), which we assume from now on.

\subsubsection{Proving Theorem \ref{the1Gen} in the case of (a), (b)}
When proving the theorem, we lose nothing when assuming that $\varsigma=1$. Indeed,  replacing $\|\cdot\|_i$ with norms which are within some factor of $\|\cdot\|_i$, the aggregation of the new norms is within the same factor of the aggregation of $\|\cdot\|_i$.
 Thus, from now on we assume that the squares
 $$\omega_i(x_i)=\|x_i\|_i^2
 $$
 of the norms $\|\cdot\|_i$
 are continuously differentiable and satisfy the relation
 \[%begin{equation}\label{eeq200}
 [\omega_i(x)+[\omega^\prime_i(x_i)]^Th_i\leq]\,\,\,\omega_i(x_i+h_i)\leq \omega_i(x)+[\nabla\omega_i(x_i)]^Th_i+\varkappa \omega_i(h_i)\;\forall x_i,h_i\in\bR^{n_i},\,1\leq i\leq K.
 \]%end{equation}
Note that $\|\nabla \omega(x)\|_{i,*}\leq 2\sqrt{\omega_i(x)}$ where $\|\cdot\|_{i,*}$ is the norm conjugate to $\|\cdot\|_i$, thus
 \begin{equation}\label{eeq122}
 |[\nabla\omega_i(x_i)]^Th_i|\leq 2\sqrt{\omega_i(x_i)\omega_i(h_i)}.
 \end{equation}

Let  $p$ be positive integer, %$\cT=\{t:0\leq t,\,\theta(t)\leq1\}$,
 and let
\begin{align*}
\overline{\cT}&=\{t\geq0:\theta(t)=1\},\\
f(x,t)&=\sum_{i=1}^K \omega_i^{p+1}(x_i)/t_i^p:\bR^{n_1+...+n_K}\times\bR^K_+\to\bR\cup\{+\infty\}\\
F(x,t)&=f^{{1\over p+1}}(x,t)\\
\phi(x)&=\min_{t\in\cT} f(x,t)\\
\Phi(x)&=\min_{t\in\cT}F(x,t)=\phi^{{1\over p+1}}(x)\end{align*}
(here, by convention, $a/0$ is 0 when $a=0$ and is $+\infty$ when $a>0$, making $f(x,t)$ lower semicontinuous on $\bR^{n_1+...+n_K}\times\bR^{K_+}$).
\\{\bf 0$^o$.} Let $q=p/(p+1)$. The function $s^2/\tau^q$ is convex in $s,\tau$ on $\{[s;\tau]:\tau\geq0\}$, and $F(x,t)=\|[\omega_1(x_1)/t_1^q;...;\omega_K(x_K)/t_K^q]\|_{p+1}$, that is, $F(x,t)$ is convex in $[x;t]\in\bR^{n_1+...+n_K}\times\bR^K_+$.
Consequently, $\phi$ and $\Phi$ are real-valued convex functions on $\bR^{n_1+...+n_K}$. Besides this, the convex and positive outside of the origin function $\Phi$ is homogeneous of degree 2 and as such is the square of a norm. Our goal in what follows is to demonstrate that with properly selected $p$ this norm can be taken as regular approximation of $\|[x_1;...;x_K]\|$  announced in Theorem \ref{the1Gen}.
\\
\textbf{1$^o$.} For evident reasons, minimizers of $f(x,t)$ over $t\in\cT$ do exist for every $x$.
We claim that given $x\neq0$, the minimizer $t(x)$ of $f(x,t)$ over $t\in\cT$ (or, which clearly is the same, over $t\in\overline{\cT}$) is unique and is continuous in $x\neq0$. Indeed, assuming w.l.o.g. that the nonzero blocks in $x$ are $x_1,...,x_m$, if $t$ and $t'$ are two distinct minimizers of $f(x,t)$
over $t\in\overline{\cT}$, then $t_i=t^\prime_i=0$ when $i>m$ and the initial blocks $\overline{t}=[t_1;...;t_m]$,  $\overline{t}^\prime=[t_1^\prime;...;t_m^\prime]$ of $t$ and $t'$  are two distinct minimizers of the function $\sum_{i=1}^m\omega_i^{p+1}(x_i)/\tau_i^p$ over $\tau$ running through the set $\{\tau:\in\bR^m:[\tau_1;...;\tau_m;0;...;0]\in\overline{\cT}\}$, which is impossible, since on the latter set the function in question is strongly convex.
To prove that $t(x)$ is continuous on the set $x\neq0$, let $\bar{x}\neq0$ and let $x^s\to \bar{x}$ as $s\to\infty$. We should lead to contradiction the assumption that the vectors $t^s=t(x^s)$ do not converge to $t(\bar{x})$. Assuming that the latter is the case and taking into account that $\overline{\cT}$ is compact, we can assume without
loss of generality that $t^s$ converge to some $t\neq t(\bar{x})$, $t\in \overline{\cT}$,  as $s\to\infty$. Since $f(x,t)$ is lower semicontinuous
on $\bR^{n_1+...+n_K}\times\bR^K_+$ and $\phi(x)$ is convex real valued and thus continuous on $\bR^{n_1+...+n_K}$, we have $\phi(x^s)=f(x^s,t^s)\to\phi(\bar{x})$
as $s\to\infty$ and therefore $f(\bar{x},t)\leq \lim_{s\to\infty} f(x^s,t^s)=\phi(\bar{x})$, implying $f(\bar{x},t)\leq\phi(\bar{x})=\min\limits_{\tau\in\overline{\cT}}f(\bar{x},\tau)$.
Thus, $t$ and $t(\bar{x})$ are two distinct minimizes of $f(\bar{x},\tau)$ over $\tau\in\overline{\cT}$, which, as we already know,  is impossible when $\bar{x}\neq0$.
\\
\textbf{2$^o$.} Observe that when $x\neq0$,  $t=t(x)$  if and only if there exists  $\lambda(x)\geq0$ such that $t$  is the unique solution to the system of equations
%\begin{subequations}
\begin{align}
t&\geq0, \;\;\theta(t)=1,\nonumber\\ %\label{eeq0.a}
x_i&=0\Rightarrow t_i=0, \nonumber \\ %\label{eeq0.b}\\
x_i&\neq0\Rightarrow p\omega_i^{p+1}(x_i)/t_i^{p+1}=\lambda(x)\underbrace{[\nabla\theta(t)]_i}_{\theta_i(x)}\label{eeq0.c}
\end{align}
%\end{subequations}
in variable $t$.
Consequently, when $x\neq0$, we have $\lambda(x)>0$. Setting $I(x)=\{i:x_i\neq0\}$, for $i\in I(x)$ we have  $\theta_i(x)>0$ and
$
t_i(x)=\big[\lambda(x)\theta_i(x)/p\big]^{-{1\over p+1}}\omega_i(x_i)
$, whence
$$
1=\theta(t(x))=\sum_it_i(x)\big[\nabla\theta(t(x))\big]_i=[\lambda(x)/p]^{-{1\over p+1}}\sum_i\omega_i(x_i)\theta_i^{{p\over p+1}}(x).
$$
Thus, when $x\neq0$, we have
\begin{subequations}
\begin{align}
\lambda(x)&=p\left[\sum_i\omega_i(x_i)\theta_i^{{p\over p+1}}(x)\right]^{p+1}>0\nonumber \\%\label{eq1.a}\\
t_i(x)&=\left\{\begin{array}{ll}{\omega_i(x_i)\theta_i^{-{1\over p+1}}(x)\over
\sum_j\omega_j(x_j)\theta_j^{{p\over p+1}}(x)},&i\in I(x)=\{i:x_i\neq0\}\\
0,&i\not\in I(x)\\
\end{array}\right.\label{eq1.b}\\
\phi(x)&=\left[\sum_i\omega_i(x_i)\theta_i^{{p\over p+1}}(x)\right]^{p+1}.\label{eq1.c}
\end{align}
\end{subequations}
\\\textbf{3$^0$.} Let $x\neq0$. For $[x;t]\in\Dom f$ let
\begin{equation}\label{eeq5}
\begin{array}{rcl}
[f^\prime_x(x,t)]_i&=&\left\{\begin{array}{ll}0,&x_i=0,\\
(p+1)\omega_i^p(x_i)\nabla\omega_i(x_i)/t_i^p(x),&x_i\neq0,\cr\end{array}\right.\\
{[f^\prime_t(x,t)]_i}&=&\left\{\begin{array}{ll}-p\omega_i^{p+1}(x_i)/t_i^{p+1}(x),&x_i\neq0,\cr
0,&x_i=0,\cr\end{array}\right.\\
f^\prime(x,t)&=&[f^\prime_x(x,t);f^\prime_t(x,t)],
\end{array}
\end{equation}
so that $f^\prime(x,t)\in\partial f(x,t)$ for all $[x;t]\in\Dom f$; let also $\phi^\prime(0)=0$.  We claim that $\phi^\prime(x);=f^\prime_x(x,t(x))$ is continuous vector field which is the subgradient (and due to continuity -- gradient) vector field of $\phi(x)$ on $\bR^n$. Indeed, by \rf{eq1.c} for $x\neq0$ we have
\[%begin{equation}\label{eeq6}
[\phi^\prime(x)]_i=\left\{\begin{array}{ll}(p+1)\theta_i^{{p\over p+1}}(x)\left[\sum_j\omega_i(x_i)\theta_j^{{p\over p+1}}(x)\right]^p\nabla\omega_i(x_i),&x_i\neq0\\
0,&x_i=0\\
\end{array}\right.,
\]%\end{equation}
and $\theta_i(x)=[\theta'(t(x))]$ is continuous in $x$ along with $t(x)$; continuity of $\phi^\prime(x)$ is evident. To prove that $\phi^\prime(x)$ is the gradient field of $\phi(x)$, due to continuity of the field it suffices to verify that
$\phi^\prime(x)\in\partial\phi(x)$ whenever $x$ has no zero blocks. For such an $x$, setting $t=t(x)$,  we obtain for $t'\in\overline{\cT}$
\begin{align*}
f(x',t')+\lambda(x)(\theta(t')-1) &\geq \underbrace{[f^\prime_x(x,t)]^T(x'-x)}_{[\phi'(x)]^T(x'-x)}\\
&+\underbrace{[f^\prime_t(x,t)]^T(t'-t)+\lambda(x)[\theta'(t)]^T(t'-t)}_{=0 \hbox{\tiny\ by \rf{eeq0.c}}}+f(x,t)+\lambda(x)(\theta(t)-1)\\
\end{align*}
whence
\[
f(x',t')\geq \phi(x)+[\phi'(x)]^T(x'-x)\]
and
\[\phi(x')=\min_{t'\in\overline{\cT}}f(x',t')\geq \phi(x)+[\phi'(x)]^T(x'-x),
\]as claimed.
\\\textbf{4$^o$.} We need the following
\begin{lemma}\label{lem1} Let $G$ be an open convex domain in $\bR^n$, $g(\cdot):G\to\bR$ be a continuously differentiable convex function, $\|\cdot\|$ be a norm on $\bR^n$ and $\gamma$ be a positive real. Assume that for every
$x\in G$ there exists $r=r_x>0$ such that the ball $V_x=\{u:\|u-x\|<r\}$ is contained in $G$ and for every $u,v\in V_x$ one has
$$
g(v)\leq g(u)+[g'(u)]^T(v-u)+\tfrac{\gamma}{2}\|u-v\|^2.
$$
Then $g'(u)$ is Lipschitz continuous with constant $\gamma$ on $G$ w.r.t. $\|\cdot\|$:
\begin{equation}\label{weconc}
\forall y,z\in G: \|g'(y)-g'(z)\|_*\leq\gamma \|y-z\|,
\end{equation}
where $\|\cdot\|_*$ is the norm conjugate to $\|\cdot\|$.
\end{lemma}
\textbf{Proof.} Let us fix $x\in G$. Applying Lemma \ref{simplem} to $V=V_x$ and $f$---restriction of $g$ on $V$, we conclude that (\ref{weconc}) holds true when $y,z\in V_x$; since $x\in G$ is arbitrary, (\ref{weconc}) holds true for all $y,z\in G$. \qed
Lemma \ref{lem1} admits the following immediate

\begin{corollary}\label{cor1} Let $g(x)$ be a continuously differentiable convex function on $\bR^n$, and let $\bR^n_o$ be the set of all vectors from $\bR^n$ with all entries different from 0. Let $\|\cdot\|$ be a norm on $\bR^n$ and $\gamma>0$ be a real. Assume that for every $x\in \bR^n_o$ there exists $r_x>0$ such that for the centered at $x$ $\|\cdot\|$-ball $V_x$ of radius $r_x$, $V_x\subset\bR^n_o$, it holds
$$
\forall u,v\in V_x: g(v)\leq g(u)+[g'(u)]^T(v-u)+\tfrac{\gamma}{2}\|u-v\|^2.
$$
Then
\begin{equation}\label{eq777}
\forall x,y\in\bR^n: \|g'(x)-g'(y)\|_*\leq \gamma\|x-y\|.
\end{equation}
\end{corollary}
\def\sign{\mathrm{sign}}
Indeed, since $g$ is continuously differentiable and $\bR^n_o$ is dense in $\bR^n$, it suffices to verify (\ref{eq777}) for $x,y\in\bR^n_o$. In this case, setting $x_t=x+t(y-x)$, $0\leq t\leq1$, we can find points $t_1<t_2<...<t_k$ in $(0,1)$ such that setting $t_0=0$, $t_{k+1}=1$, the intervals $\Delta_s=\{x_t:t_s<t<t_{s+1}\}$, $s=0,1,...,k$,
belong to $\bR^n_o$, that is, $\sign([x_t]_i)=\epsilon_i(s)\in\{-1,1\}$, $t_s<t<t_{s+1}$.
Then $G_s=\{x\in\bR^n:\sign(x_i)=\epsilon_i(s),i\leq n\}$ is an open convex domain containing $\Delta_s$.
By Lemma \ref{lem1}, when $t_s<t'<t''<t_{s+1}$ we have
\[\|g'(x_{t''})-g'(x_{t'})\|_*\leq \gamma\|x_{t''}-x_{t'}\|=\gamma[t''-t']\|y-x\|.
\] Passing to limit as $t'\to t_s+0$ and $t''\to t_{s+1}-0$ and taking into account that $g$ is C$^1$ on $\bR^n$, we get \[\|g'(x_{t_{s+1}})-g'(x_{t_s})\|_*\leq\gamma [t_{s+1}-t_s]\|y-x\|,\]
whence $\|g'(x_{t_{k+1}})-g'(x_{t_\aic{1}{0}})\|_*\leq\gamma\|y-x\|$. \qed
\par\noindent\textbf{5$^0$.}
\aic{Observe that {\color{red} with  $\psi(s)=s^{{1\over p+1}}:\bR\to \bR_{+}\cup\{-\infty\}$, $\bar{s}>0$ and $\delta\in\bR$ ?????}
we have
$$
\begin{array}{ll}
\psi(\bar{s}+\delta)\leq \psi(\bar{s})+{1\over p+1}\bar{s}^{-{p\over p+1}}\delta
=\psi(\bar{s})+{1\over p+1}{\delta\over \psi^p(\bar{s})}.\\
\end{array}
$$}{}
Let now $x\in\bR^{n_1+...+n_K}$ be a vector with all entries different from 0 and $h\in\bR^{n_1+...+n_K}$ be a vector with \be
\omega^{1/2}_i(h_i)<\omega^{1/2}_i(x_i)/[\aic{2(p+1)}{(2p+2)}^2(1+\varkappa)],\;\; i\leq K.
\ee{eqhi<}
One has
$$
\omega_i(x_i)+[\nabla\omega_i(x_i)]^Th_i\leq \omega_i(x_i+h_i)\leq \omega_i(x_i)+[\nabla\omega_i(x_i)]^Th_i+\varkappa\omega_i(h_i).
$$
Thus, when setting
$$
\Delta_i=[\omega_i(x_i+h_i)-\omega_i(x_i)]/\omega_i(x_i)
$$
we have
\[
\Delta_i=[\nabla\omega_i(x_i)^Th_i]/\omega_i(x_i)+\delta_i,\quad |\delta_i|\leq\varkappa\omega_i(h_i)/\omega_i(x_i),
\]
and taking into account (\ref{eeq122}) and \rf{eqhi<},
\[ |\Delta_i|\leq 2\sqrt{\omega_i(h_i)/\omega_i(x_i)}+\varkappa\omega_i(h_i)/\omega_i(x_i)\leq 3\sqrt{\omega_i(h_i)/\omega_i(x_i)}\leq 3/(2p+2)^2.
\]
We have
\begin{align*}
&\phi(x)+[\phi'(x)]^Th\leq \phi(x+h)\leq \sum_i\left[\sum_i\omega_i^{p+1}(x_i+h_i)\right]/t_i^p(x)
= \sum_i\omega_i^{p+1}(x_i)(1+\aic{\Delta}{\Delta_i})^{p+1}\\
&\leq \sum_i\omega_i^{p+1}(x_i)\left[1+(p+1)\Delta_i+{p(p+1)\over 2}\Delta_i^2\left[1+\sum_{j=1}^{p-1}\aic{|p\Delta_i|^j}{|(p-1)\Delta_i|^j\over j!}\right]\right]t_i^{-p}(x)\\
&\leq \sum_i\omega_i^{p+1}(x_i)\left[1+(p+1)\Delta_i+\aic{2p(p+1)}{{p(p+1)\over 2}\exp\left\{(p-1)|\Delta_i|\right\}}\Delta_i^2\right]t_i^{-p}(x)\\
&\leq \sum_i\omega_i^{p+1}\bigg[1+(p+1)\left[{\nabla\omega_i(x_i)^Th_i\over \omega_i(x_i)}+\varkappa{\omega_i(h_i)\over\omega_i(x_i)}\right]
+\aic{9(p+1)(2p+1)}{5p(p+1)}{\omega_i(h_i)\over\omega_i(x_i)}\bigg]t_i^{-p}(x)\\
&=\phi(x)+[\phi'(x)]^Th+\aic{9(p+1)(2p+1+\varkappa)}{(p+1)(5p+\varkappa)}\sum_i\omega_i(h_i)\omega_i^p(x_i)t_i^{-p}(x)\\
&=\phi(x)+\underbrace{[\phi'(x)]^Th+\aic{9(p+1)(2p+1+\varkappa)}{(p+1)(5p+\varkappa)}\left[{\sum}_i\omega_i(h_i)\theta_i^{{p\over p+1}}(x)\right]\left[{\sum}_i\omega_i(x_i)\theta_i^{{p\over p+1}}(x)\right]^p}_{\delta}
\end{align*}
where the concluding equality is due to (\ref{eq1.b}).

\hide{Note that
$$
\begin{array}{rcl}
\bar{s}&:=&\phi(x)=\left[\sum_i\omega_i(x_i)^2\theta_i^{{p\over p+1}}(x)\right]^{p+1}\\
\delta&=&[\phi'(x)]^Th+9(p+1)(2p+1)\left[{\sum}_i\omega_i(h_i)^2\theta_i^{{p\over p+1}}(x)\right]\left[{\sum}_i\omega_i(x_i)^2\theta_i^{{p\over p+1}}(x)\right]^p\\
&=&\left[\sum_i2(p+1)[x_i^Th_i+{9(2p+1)\over2}\omega_i(h_i)^2]\theta_i^{{p\over p+1}}(x)\right]\left[\sum_j\|x_j\|_2^2\theta_j^{{p\over p+1}}(x)\right]^p\\
&\geq&-{1\over2p+2}\bar{s} \hbox{\ [since $\omega_i(h_i)\leq \omega_i(x_i)/(2p+2)^2$]}\\
\end{array}
$$
implying that
$\phi(x+h)\leq \phi(x)+\delta$, whence}{
Note that for \[\psi(s)=\left\{\begin{array}{ll}s^{{1\over p+1}}& s\geq 0,\\-\infty &s<0,\end{array}\right.
\] $\bar{s}>0$ and $\delta\in\bR$ one has
\[\psi(\bar{s}+\delta)\leq \psi(\bar{s})+{1\over p+1}\bar{s}^{-{p\over p+1}}\delta
=\psi(\bar{s})+{1\over p+1}{\delta\over \psi^p(\bar{s})}.
\]
Hence,}
\begin{align}
&\Phi(x)+[\Phi'(x)]^Th\leq \Phi(x+h)=\psi(\phi(x+h))
\leq \psi(\phi(x)+\delta)\leq\psi(\phi(x))+{1\over p+1}\delta/\psi^p(\phi(x))\nn
&=\Phi(x)+{1\over p+1}[\phi'(x)]^Th/\psi^p(\phi(x))\nn
&\quad\quad+\aic{9(2p+1+\varkappa)}{(5p+\varkappa)}\left[{\sum}_i\omega_i(h_i)\theta_i^{{p\over p+1}}(x)\right]\left[{\sum}_i\omega_i(x_i)\theta_i^{{p\over p+1}}(x)\right]^p\psi^{-p}(\phi(x))\nn
&=\Phi(x)+[\Phi'(x)]^Th+\aic{9(2p+1+\varkappa)}{(5p+\varkappa)}\left[{\sum}_i\omega_i(h_i)\theta_i^{{p\over p+1}}(x)\right]\left[{\sum}_i\omega_i(x_i)\theta_i^{{p\over p+1}}(x)\right]^p\psi^{-p}(\phi(x))\nn
&=\Phi(x)+[\Phi'(x)]^Th+\aic{9(2p+1+\varkappa)}{(5p+\varkappa)}\left[{\sum}_i\omega_i(h_i)\theta_i^{{p\over p+1}}(x)\right],\label{eq2}
\end{align}
due to $\psi^p(\phi(x))=\phi^{{p\over p+1}}(x)=\left[\sum_i\omega_i(x_i)\theta_i^{{p\over p+1}}(x)\right]^p$ by (\ref{eq1.c}).
\par\noindent\textbf{6$^0$.} Now, for $x\in\bR^{n_1+...+n_K}$, let us set $\omega[x]=[\omega_1(x);...;\omega_K(x)]$ and $\|x\|=\theta^{1/2}(\omega[x])$. Clearly, $\|\cdot\|$ indeed is a norm, and the set $\{x:\|x\|=1\}$ is nothing but the set $\{x:\omega[x]\in\overline{\cT}\}$. \par
Given $\beta\geq1$, we say that two nonnegative reals $a,b$ are within factor $\beta$, if $\beta^{-1}a\leq b\leq\beta a$, or, which is the same, $\beta^{-1}b\leq a\leq \beta b$. Let us say that two block vectors $u,v\in\bR^{n_1+...+n_K}$ are within factor $\beta$, if $\sqrt{\omega_i(u_i)}$ and $\sqrt{\omega_i(v_i)}$ are within factor $\beta$ for every $i$.
\par Let us make the following observation
\begin{quote}
{\sl Let $x\in\bR^{n_1+...+n_K}$ be a vector with nonzero blocks such that $\|x\|=1$. Then there exists vector $z$ with nonzero blocks such that $\phi(z)=1/p$ and $z$ is within factor
$$
\kappa=[pK]^{{1\over 2(p+1)}}
$$
from $x$. %; recall that the entries in  $\nabla \theta(t)$, for every $t\in\overline{\cT}$, are within factor $K$ from 1.
Vice versa, if $z$ is a vector with nonzero blocks such that $\phi(z)=1/p$, then there exists vector $x$ with nonzero blocks such that $\|x\|=1$ and $x$ is within factor $\kappa$ from $z$.
}\end{quote}
Indeed, let $x$ be vector with nonzero blocks such that $\|x\|=1$, and let $t=\omega[x]$, so that $t\in\overline{\cT}$ and $t>0$. Let now $z$ be with nonzero blocks $z_i\in\bR^{n_i}$ such that $\sqrt{\omega_i(z_i)}=\left[p^{-1}t_i^{p+1}[\nabla \theta(t)]_i\right]^{{1\over 2(p+1)}}$, so that
\[\omega_i(z_i)=t_i[p^{-1}[\nabla \theta(t)]_i]^{{1\over p+1}}=\omega_i(x_i)[p^{-1}[\nabla \theta(t)]_i]^{{1\over p+1}}.
\] Note that $z$ and $x$ are within the factor $\kappa$. We have $t\in\overline{\cT}$ and $p\omega_i^{p+1}(z_i)/t_i^{p+1}=[\nabla \theta(t)]_i$, implying that $t=t(z)$ by (\ref{eeq0.c}). Consequently,
$$
\phi(z)=\sum_i\omega_i^{p+1}(z_i)/t_i^p=p^{-1}\sum_it_i[\nabla \theta(t)]_i=p^{-1}\theta(t)=p^{-1},
$$
as  claimed. Vice versa, let $z$ be a vector with nonzero blocks such that $\phi(z)=1/p$, so that for $t=t(z)\in\overline{\cT}$ and  $\lambda=\lambda(z)\geq0$ we have
\begin{align*}
\omega_i^{p+1}(z_i)&=(\lambda/p)t_i^{p+1}[\nabla \theta(t)]_i,\\
 1/p=\phi(z)&=\sum \omega_i^{p+1}(z_i)/t_i^p=(\lambda/p)\sum_it_i[\nabla \theta(t)]_i=(\lambda/p)\theta(t)=\lambda/p.
\end{align*}
We see that $\lambda=1$, and therefore $\omega_i(z_i)$ are within factor $\kappa^2$ from $t_i$. Selecting $x_i\in\bR^{n_i}$ such that $\omega_i(x_i)=t_i$, we get vector with nonzero blocks such that $\|x\|=\theta(t)=1$ and $z$ is within factor $\kappa$ from $x$. \qed

Next, $\Phi(x)$ is convex positive outside of the origin homogeneous, of degree 2, function and thus  $\mfn(x):=\sqrt{\Phi(x)}=\phi^{{1\over 2(p+1)}}(x)$ is a norm.
By construction $\mfn(\cdot)$ is an ``absolute block-norm,'' meaning that $\mfn(x)$ depends solely and monotonically on the   norms $\|x_i\|_i=\sqrt{\omega_i(x_i)}$  of the blocks $x_i$.
As an immediate corollary, if $u,v$ are two vectors with nonzero blocks which are within some factor from each other,
then $\mfn(u)$ and $\mfn(v)$ are within the same factor from each other as well. Given vector $x$ with nonzero blocks and such that $\|x\|=1$,
 there exists vector $z$ within factor $\kappa$ from $x$ with $\phi(z)=1/p$ and thus with $\mfn(z)=\pi:=p^{-{1\over 2(p+1)}}$, so that $\mfn(x)$ is within factor $\kappa$ from $\pi$. Thus,
$\|\cdot\|$ is within factor $\overline{\kappa}=\kappa/\pi$ from $\mfn(\cdot)$.
\\
\textbf{7$^o$.} Our next observation is as follows:
\begin{quote}{\sl
Let $t\in\overline{\cT}$ and $h\in\bR^n$. Then
$$
\sum_i\omega_i(h_i)[\nabla \theta(t)]_i^{{p\over p+1}}\leq K^{{1\over p+1}} \|h\|^2.
$$
}\end{quote}
Indeed, by (\ref{eeq334}) we have $\hide{{\color{red}\theta_i^{{p\over p+1}}(t)}}{[\nabla \theta(t)]_i^{{p\over p+1}}} \leq K^{{1\over p+1}}[\nabla \theta(t)]_i$, whence
\begin{align*}
\sum_i\omega_i(h_i)[\nabla \theta(t)]_i^{{p\over p+1}}&\leq K^{{1\over p+1}}[\nabla \theta(t)]^T\omega[h]\leq K^{{1\over p+1}}\left[\theta(t+\omega[h])-\theta(t)\right]\\
&\leq K^{{1\over p+1}}[\theta(t)+\theta(\omega[h])-\theta(t)]
=K^{{1\over p+1}}\|h\|^2.
\end{align*}
As a consequence, we have
$$
\sum_i\omega_i(h_i)[\nabla \theta(t)]_i^{{p\over p+1}}\leq K^{{1\over p+1}}(\kappa/\pi)^2\mfn^2(h)
$$
Finally, invoking (\ref{eq2}), we conclude that for every vector $x$ with nonzero blocks and all $h$ such that
\[\omega_i^{1/2}(h_i)\leq \omega_i^{1/2}(x_i)/[(2p+2)^2(1+\varkappa)],\;\; i\leq K,
\] one has
\[
\mfn^2(x+h)\leq\mfn^2(x)+[\nabla \mfn^2(x)]^Th+\tfrac{\gamma}{2}\mfn^2(h),\;\;\gamma=\aic{18K^{{1\over p+1}}(\kappa/\pi)^2(2p+1+\varkappa)}{2K^{{1\over p+1}}(\kappa/\pi)^2(5p+\varkappa)}.
\]
We arrive at the following conclusion:
\begin{quote}
{\sl In \aic{the Special case,}{ the situation in question,} for every positive integer $p$ there exists a norm $\mfn(\cdot)$ on $\bR^{n_1+...+n_K}$ such that
\begin{itemize}
\item  $\Phi(x)=\mfn^2(x)$ is continuously differentiable with Lipschitz continuous gradient and satisfies
$$
\forall x,h\in\bR^n: \Phi(x+h)\leq\Phi(x)+[\nabla\Phi(x)]^Th+\aic{9[pK]^{{3\over 2(p+1)}}(2p+1+\varkappa)}{[pK]^{{2\over p+1}}(5p+\varkappa)}\Phi(h)
$$
\item $\mfn(\cdot)$ is within the factor $\kappa/\pi=K^{{1\over 2(p+1)}}p^{{1\over p+1}}$ from the norm $\|x\|=\sqrt{\theta(\omega[x^2])}$.
\end{itemize}}
\end{quote}
Setting  $p=\lceil\ln(K+1)\rceil$\footnote{Here $\lceil a\rceil$ stands for the upper integer part of $a$---the smallest integer greater or equal to $a$.} we complete the proof of Theorem \ref{the1Gen}.\qed

\subsection{Proof of \aic{Theorem}{Proposition} \ref{thefact}}
As is immediately seen, we lose nothing when assuming that the mapping $z\to Pz$ in (\ref{eeq100}) is just the natural projection $\bR^n=\bR^{p+q}\ni[u;y]\mapsto u\in\bR^p$. In this case, factor-norm induced by this projection and the norm $\|\cdot\|$ on $\bR^{p+q}$:
$$
\|u\|'=\min_y\|[u;y]\|.
$$
Observe that when two norms on $\bR^{p+q}$ are within some factor from each other, their factor-norms induced by $P$ are within the same factor as well. Thus, it suffices to prove Theorem \ref{thefact} in the case when the function $\Phi(x)=\|x\|^2$ if continuously differentiable and such that
$$
\|\nabla\Phi(x)-\nabla\Phi(x')\|_*\leq 2\varkappa\|x-x'\|\,\,\forall x,x';
$$
where, as always, $\|\cdot\|_*$ is the norm conjugate to $\|\cdot\|$.
In this situation the statement of \aic{Theorem}{Proposition} \ref{thefact} is readily given by the following
\begin{lemma}\label{lemfact} Let $\Psi(u,y):\bR^{p+q}$ be a convex function, $\|\cdot\|$ be a norm on $\bR^{p+q}$, and let $\|\cdot\|'$ be the factor-norm:
$$
\|h\|'=\min_{d\in\bR^q}\|[h;d]\|.
$$
induced by the natural projection of $\bR^{p+q}$ onto $\bR^p$.
 Assume that $\Psi$ is smooth with Lipschitz continuous with constant $L$ w.r.t. $\|\cdot\|$ gradient, i.e.,
\[\|\nabla\Psi(z)-\nabla\Psi(z')\|_*\leq L\|z-z'\|\quad\forall z,z',
\]
and let the set
$$
Y(u)=\Argmin_{y\in\bR^q}\Psi(u,y)
$$ be nonempty for every $u\in\bR^p$. Then the function
$$
\overline{\Psi}(u)=\min_y\Psi(u,y)
$$
 on $\bR^p$ is convex and smooth with Lipschitz continuous with constant $L$ w.r.t. $\|\cdot\|'$ gradient.
\end{lemma}
{\bf Proof of the lemma.} Convexity of $\overline{\Psi}$ is evident. Let $u\in \bR^p$, and let $y\in Y(u)$. Then for all $(u',y')$,
$$\Psi(u',y')\geq \Psi(u,y)+[\Psi^\prime_u(u,y)]^T(u'-u) +[\Psi^\prime_y(u,y)]^T(y'-y)=\Psi(u,y)+[\Psi^\prime_u(u,y)]^T(u'-u),
$$
that is $\Psi(u',y')\geq\overline{\Psi}(u)+[\Psi^\prime_u(u,y)]^T(u'-u)$ for all $u',y'$; hence,
\begin{equation}\label{eeq34}
\overline{\Psi}(u')\geq  \overline{\Psi}(u)+[\Psi^\prime_u(u,y)]^T(u'-u)\,\forall u'\in\bR^p.
\end{equation}
Denoting for brevity $g=\Psi^\prime_u(u,y)$, for every $h\in\bR^p$, selecting $d\in\bR^q$ such that $\|[h;d]\|=\|h\|'$,  we get
\[%begin{equation}\label{eeq55}
\overline{\Psi}(u+h)\leq \Psi(u+h,y+d)\leq \Psi(u,y)+g^Th+\tfrac{L}{2}\|[h;d]\|^2=\Psi(u,y)+g^Th+\tfrac{L}{2}[\|h\|']^2,
\]%end{equation}
implying that
\begin{equation}\label{eeq111}
\overline{\Psi}(u)+g^Th\leq \overline{\Psi}(u+h)\leq \overline{\Psi}(u)+g^Th+{L\over 2}[\|h\|']^2
\end{equation}
with the first inequality given by (\ref{eeq34}).
Consequently, $g$ is the Frechet derivative of $\overline{\Psi}$ at $u$ and is independent of how $y\in Y(u)$ is selected. Next, $\overline{\Psi}$ is convex
on $\bR^p$ and thus is locally Lipschitz continuous. We claim that $\overline{\Psi}'(\cdot)$ is continuous. Indeed, from local Lipschitz continuity of $\overline{\Psi}$
it follows that $\overline{\Psi}'(\cdot)$ is locally bounded. Thus, to prove the continuity of $\overline{\Psi}'(\cdot)$, it suffices to lead to the contradiction the assumption that for some sequence $u_i\to u$, $i\to\infty$, such that $g_i=\overline{\Psi}'(u_i)$ converge to some $e$, one has $e\neq g:=\overline{\Psi}'(u)$. But this is evident: $e$ clearly is a subgradient of $\overline{\Psi}$ at $u$, and since $\overline{\Psi}$ is differentiable, we get $e=g$, which is a contradiction.\par
Thus, $\overline{\Psi}$ is continuously differentiable convex function satisfying (\ref{eeq111}), implying by Lemma \ref{lem1} that the gradient of $\overline{\Psi}$ is Lipschitz continuous, with constant $L$, w.r.t. $\|\cdot\|'$. \qed

\subsection{Proof of Theorem \ref{theass}}
Since $\theta$ is $(\overline{\varkappa},\overline{\varsigma})$-regular and $\|\cdot\|_i$ are $({\varkappa}',{\varsigma}')$-regular, there exist $\overline{\varkappa}$-smooth norm $\vartheta(\cdot)$ and ${\varkappa}'$-smooth norms $\vartheta_i$ such that
\begin{equation}\label{eq456}
\overline{\varsigma}^{-1}\vartheta(\cdot)\leq \theta(\cdot)\leq \overline{\varsigma} \vartheta(\cdot),\,\,(\varsigma')^{-1}\vartheta_i(\cdot)\leq \|\cdot\|_i\leq\varsigma'\vartheta_i(\cdot),\,i\leq K.
\end{equation}
Since $\theta$ is an absolute norm, these relations hold true when replacing $\vartheta(\cdot)$ with the absolute norm
$$
\overline{\vartheta}(y)=2^{-K} \sum_{E\in{\cE}_K}\vartheta(Ey)
$$
where ${\cal E}_K$ is the multiplicative group of diagonal $K\times K$ matrices with diagonal entries $\pm1$, and the norm $\overline{\vartheta}$ clearly is $\overline{\varkappa}$-smooth along with $\vartheta$. Thus, we can assume, on the top of what was just stated, that $ \vartheta$ is an absolute norm.
\par
Next, the function $\Psi(x_1,...,x_K)={\theta^{1/2}(\|x_1\|_1^2,...,\|x_K\|_K^2)}$ is positively homogeneous of degree 1, and since $\theta$ is an absolute norm, function $\Psi^2(y)$ is monotone on $\bR^K_+$. Consequently, $\Psi^2$ is convex, and therefore the Lebesque sets of $\Psi$ are convex; together with positive homogeneity of degree 1 of the function this implies that $\Psi(x_1,...,x_K)$ is a norm. By the same argument, function ${\vartheta^{1/2}(\vartheta_1^2(x_1),...,\vartheta_K^2(x_K))}$ also is a norm, so that monotonicity and homogeneity of degree 1 of $\theta$ and $\vartheta$ on $\bR^K_+$ combines with (\ref{eq456}) to imply that for $\varsigma={\varsigma}'\sqrt{\overline{\varsigma}}$ and  all $[x_1;...;x_K]\in\bR^{n_1+...+n_K}$
\[\varsigma^{-1}{\vartheta^{1/2}(\vartheta_1^2(x_1),...,\vartheta_K^2(x_K))}\leq {\theta^{1/2}(\|x_1\|_1^2,...,\|x_K\|_K^2)}\leq
\varsigma{\vartheta^{1/2}(\vartheta_1^2(x_1),...,\vartheta_K^2(x_K))}.
\]
To complete the proof of the theorem, all we need is the following statement.
\begin{lemma}\label{thenew} Let $F(y)=\|y\|_o^2$ be the square of an absolute norm on $\bR^K$ such that $\nabla F$ is  Lipschitz continuous w.r.t. $\|\cdot\|_o$ with Lipschitz constant $L$. Let also $\phi_i(x_i)=\|x_i\|_i^2$ be squares of norms on $\bR^{n_i}$ such that $\nabla \phi_i$ is Lipschitz continuous with constant $\ell$ with respect to $\|\cdot\|_i$, $1\leq i\leq K$. For $x=[x_1;...;x_K]\in\bR^{n_1}\times...\times\bR^{n_k}$, we set $\phi(x)=[\phi_1(x_1);...;\phi_K(x_K)]$ and
\[
f(x)=F^{1/2}(\phi(x)).
\]
Then $f$ is the square of a norm $\|\cdot\|$ on $\bR^{n_1+...+n_K}$ with gradient $\nabla f$ which is Lipschitz continuous with constant $2L+\ell$ w.r.t.  $\|\cdot\|$.
\end{lemma}
{\bf Proof of Lemma \ref{thenew}.} The fact that $f$ is square of a norm is evident---since $F$ is square of an absolute norm the function $F(\phi(x))$ is convex. Besides this, the function is continuously differentiable, and is of homogeneity degree 4; as a result, $F^{1/4}(\phi(x))$ is continuous positively homogeneous, of degree 1, even function with convex \aic{Lebesque}{level} sets and as such is convex, and therefore is a norm. Note also that $f$ is continuously differentiable. All we need to prove is that
\[%begin{equation}\label{eqq1}
\forall (x,h): f(x+h)\leq f(x)+[f'(x)]^Th+\half[2L+\ell]f(h).
\]%end{equation}
Assume for the sake of simplicity that $F$ and $\phi_i$ are twice continuously differentiable outside of  the respective origins (we can arrive at this situation by suitable approximation, see Lemma \ref{lemtech} below).
Invoking \aic{Corollary \ref{cor1}}{Lemma \ref{simplem}}, it suffices to show that for every $x$ with nonzero $x_i$'s and every $h$ it holds
\begin{equation}\label{eqq2}
D^2f(x)[h,h]\leq [2L+\ell]f(h).
\end{equation}
Let us put
$$
\cD=\{x\in\bR^{n_1+...+n_K}:x_i\neq0,i\leq K\},\;\;\;\cD_o=\{y\in\bR^K:y\neq0\}.
$$
$F$ is twice continuously differentiable on $\cD_o$ and
\begin{align*}
0&\leq D^2F(y)[d,d]\leq LF(d)\;\;\forall (y\in\cD_o,d\in\bR^K),&(a)\\
0&\leq D^2\phi(x)[h,h]\leq \ell \phi(h)\;\;\forall (x\in\cD,h\in\bR^{n_1+...+n_K}).&(b)
\end{align*}
Besides this, taking into account that $F$ and $\phi_i$ are squares of norms continuously differentiable outside of  the origin, we have
\def\abs{\hbox{abs}}
\begin{align*}
\abs[D\phi(x)[h]]&\leq 2[\phi(x)]^{1/2}\cdot[\phi(h)]^{1/2},&(c)\\
|DF(y)[d]|&\leq 2\sqrt{F(y)}\sqrt{F(d)},&(d)\\
y\geq0&\Rightarrow DF(y)[h]\hbox{\ is monotone in $h$}&(e)
\end{align*}
where for a vector $z$, $\abs[z]$ is the vector of modulae of entries in $z$,
 for a nonnegative vector $z$ and $\alpha>0$, $[z]^\alpha$ stands for the vector comprised of $\alpha$-th degrees of the entries of $z$, and
 $u\cdot v$ is coordinatewise product of two vectors; \aic{$(d)$}{$(e)$} stems from the fact that $F$ is the square of an absolute norm.
\par
Now, let $x\in\cD$ and let $y=\phi(x)$, so that $y\in \cD_o$. We have $f(x)=F^{1/2}(\phi(x))$ with
\[Df(x)[h]={1\over 2F^{1/2}(\phi(x))}DF(\phi(x))[D\phi(x)[h]]
\] and
\begin{align}
D^2f(x)[h,h]&={1\over 2F^{1/2}(\phi(x))}\left[D^2F(\phi(x))[D\phi(x)[h],D\phi(x)[h]]+DF(\phi(x))[D^2\phi(x)[h,h]]\right]\nn
&\quad-{1\over 4F^{3/2}(\phi(x))}(DF(\phi(x))[D\phi(x)[h]])^2\nn
&\leq{1\over 2F^{1/2}(y)}\left[D^2F(y)[D\phi(x)[h],D\phi(x)[h]]+DF(y)[D^2\phi(x)[h,h]]\right].
\label{eqd2f}
\end{align}
By homogeneity, it suffices to verify (\ref{eqq2}) for $x\in\cX$ such that $F(\phi(x))=1$.
In this case we have
$$
\begin{array}{l}
Df(x)[h]={1\over 2}DF(y)[D\phi(x)[h]]\\
D^2f(x)[h,h]\leq{1\over 2}\left[D^2F(y)[D\phi(x)[h],D\phi(x)[h]]+DF(y)[D^2\phi(x)[h,h]]\right]\\
\end{array}
$$
Putting $d=D^2\phi(x)[h,h]$, we get by $(d)$,
\begin{align*}
DF(y)[D^2\phi(x)[h,h]]&\leq 2\sqrt{F(D^2\phi(x)[h,h])}\leq 2\ell \sqrt{F(\phi(h))}=2\ell f\aic{(\phi(h))}{(h)}
\end{align*}
by $(b)$ and due to homogeneity of degree 2 of $F$ and monotonicity of $F(\cdot)$ on $\bR^K_+$.
Besides this,
\begin{align*}
&D^2F(y)[D\phi(x)[h],[D\phi(x)[h]]\leq L F(D\phi(x)[h]) \hbox{\ [by $(a)$]}\\
&= LF(\abs[D\phi(x)[h]])\hbox{\ [since $F(\cdot)$ is the square of an absolute norm]}\\
&\leq LF(2[\phi(x)]^{1/2}\cdot[\phi(h)]^{1/2})\hbox{\ [by  $(c)$ and due to monotonicity of $F(\cdot)$ on $\bR^K_+$]}\\
&\leq LF(\lambda \phi(x)+\lambda^{-1}\phi(h))\,\,\forall \lambda>0\hbox{\ [since $F$ is monotone on $\bR^K_+$]}\\
%&=L[f(\lambda\phi(x)+\lambda^{-1}\phi(h))]^2\,\forall \lambda>0\\ %\hbox{\ [since $F=f^2$]}\\
&\aic{=}{\leq }L[\aic{\lambda f(\phi(x))+\lambda^{-1}f(\phi(h))}{\lambda F^{1/2}(\phi(x))+\lambda^{-1}F^{1/2}(\phi(h))}]^2\;\;\forall\lambda>0\hbox{\ [since $\aic{f}{F^{1/2}}(\cdot)$ is a norm on $\bR^K$],}
\end{align*}
whence
\begin{align*}
D^2F(y)[D\phi(x)[h],D\phi(x)[h]]&\leq L\left[\inf_{\lambda>0}[\aic{\lambda f(\phi(x))+\lambda^{-1}f(\phi(h))}{\lambda F^{1/2}(\phi(x))+\lambda^{-1}F^{1/2}(\phi(h))}]\right]^2\\
&=
4L\aic{f(\phi(x))f(\phi(h))}{F^{1/2}(\phi(x)) F^{1/2}(\phi(h))}=4Lf\aic{(\phi(h))}{(h).}
\end{align*}
The bottom line is that when $x\in\cD$ and $F(\phi(x))=1$, we have by \rf{eqd2f}
\begin{align*}
D^2f(x)[h,h]&\leq
\half\left[D^2F(y)[D\phi(x)[h],D\phi(x)[h]]+DF(y)[D^2\phi(x)[h,h]]\right] \\
&\leq \half[4Lf\aic{(\phi(h))}{(h)}+2\ell f\aic{(\phi(h)}{(h)}]=[2L+\ell]f\aic{(\phi(h))}{(h)},
\end{align*}
as required in (\ref{eqq2}). \qed
\subsection{Proof of Theorem \ref{theellispec}}
{\em (i):} By Theorem \ref{the1Gen} as applied to the (1,1)-regular norms $\|\cdot\|_i=\|\cdot\|_2$ on $\bR^n$, the norm \[\|[x_1;...;x_K]\|:=\theta^{1/2}(\|x_1\|_2^2,...,\|x_K\|_2^2)
\] on $\bR^{nK}$ is $(O(1)\ln(K+1),O(1))$-regular. With $\cX$ given by (\ref{eeq1}), the norm
\[\|x\|_+:=\theta^{1/2}(\|T_1^{1/2}x\|_2^2,...,\|T_k^{1/2}x\|_2^2)
\] is the superposition of $\|\cdot\|$ and linear embedding $x\mapsto[T_1^{1/2}x;... ;T_k^{1/2}x]$ and therefore is regular with the same parameters as $\|\cdot\|$, see Section \ref{regnorms}. Finally,  $\|\cdot\|_{\overline{\cX}}$ is the factor-norm of $\|\cdot\|_+$, and passing from a norm to its factor-norm preserves regularity parameters by Proposition \ref{thefact}.
\par\noindent
{\em (ii):} The norms $\|\cdot\|_i=\|\cdot\|_{2,2}$ on $\bS^{d_i}$ are $(O(1)\ln(d_i+1),O(1))$-regular, see ``basic facts'' at the end of in Section \ref{regnorms}.
By Theorem \ref{the1Gen}) as applied to these norms, the norm
\[\|[y_1;...;y_K]\|:=\theta^{1/2}(\|y_1\|_{2,2}^2,...,\|y_K\|_{2,2}^2)
\] on $\bS^{d_1}\times...\times\bS^{d_K}$ is $(O(1)[\ln(K+1)+\max_{i\leq K}\ln(d_i)],O(1))$-regular.
Similarly to the case of $(i)$, with $\cX$ given by \aic{(\ref{eeq100})}{\rf{eeq10043}}, the norm $\|\cdot\|_\aic{\cX}{+}$ \aic{}{which}
is the superposition of $\|\cdot\|$ and linear embedding $x\mapsto(S_1(x)_,...,S_K(x))$ \aic{and  therefore}{} is regular with the same parameters as $\|\cdot\|$, and these parameters are inherited by the factor-norm $\|\cdot\|_{\overline{\cX}}$. \qed
\subsection{Technical lemmas}
\begin{lemma}\label{simplem} Let $\|\cdot\|$ be a norm on $\bR^n$ with conjugate norm $\|\cdot\|_*$. Let also $\alpha\geq0$, let $V$ be a nonempty open convex domain in $\bR^n$ and $f:\,V\to\bR$ be a continuously differentiable convex function.  The following 3 properties of $f$ are equivalent:
\begin{enumerate}
\item[{\rm (1)}] $h^T[\nabla f(x+h)-\nabla f(x)]\leq \alpha\|h\|^2$ $\forall x,x+h\in V$;
\item[{\rm (2)}] $\|\nabla f(x+h)-\nabla f(x)\|_*\leq \alpha \|h\|$ $\forall x,x+h\in V$;
\item[{\rm (3)}] $f(x+h)\leq f(x)+h^T\nabla f (x)+{\alpha\over 2}\|h\|^2$ $\forall x,x+h\in V$.
\end{enumerate}
\end{lemma}
{\bf Proof.} Let $V_o$ be an open convex domain such that for some $\epsilon[V_o]>0$ the $\|\cdot\|_2$-norm $\epsilon$-neighbourhood of $V_o$ is contained in $V$. It suffices to demonstrate that the properties (1$_o$)--(3$_o$) obtained from respective properties (1)--(3) by replacing the requirement  $x,x+h\in V$ with $x,x+h\in V_o$ are equivalent to each other, whatever be $V_o$ of the indicated type. Let us fix such a $V_o$ and note that when $f$ satisfies property ($i$), $i=1,2,3$, then  property  ($i_o$) holds true for the restrictions onto $V_o$ of all shifts $f(\cdot-u)$ of $f(\cdot)$ with $\|u\|_2<\epsilon[V_o]$, same as for convex combinations of these shifts. Selecting a C$^\infty$ nonnegative kernel $\delta(\cdot)$ with unit integral which vanishes outside of the $\epsilon[V_o]/2$-neighbourhood of the origin and setting
$$
f_\delta(x)=\int f(x-u)\delta(u)du:V_o\to\bR,
$$
we obtain a C$^\infty$ convex function on $V_o$ which possesses properties ($i_o$), $i=1,2,3$,  provided that $f$ possesses property $(i)$. When $\delta\to+0$, functions $f_\delta$ converge to $f$ uniformly on compact subsets of $V_o$ along with their gradients. Therefore, when $f_\delta$ possesses one of the properties $(i_o)$, so does $f$, and vice versa. The bottom line is that all we need is to show the equivalence $(1_o)$ -- $(3_o)$ in the case when $f$ is a convex C$^\infty$ function on $V_o$. This is immediate: when $f$ is convex and smooth, every one of the three properties in question is equivalent to
\[
u^T\nabla^2f(x)v\leq \alpha \|u\|\|v\|\;\;\forall (x\in V_o,u,v\in\bR^n).\tag{$*_o$}
\]
Indeed,
\par
$(1_o)\Leftrightarrow (*_o)$: When $f$ obeys $(1_o)$, $x\in V_o$ and $h\in\bR^n$, we clearly have
$$
h^T\nabla^2f(x)h=\lim_{t\to +0}t^{-1}h^T[\nabla f(x+th)-\nabla f(x)]=\lim_{t\to+0}t^{-2}\underbrace{[th]^T[\nabla f(x+th)-\nabla f(x)]}_{\leq\alpha\|th\|^2}\leq \alpha\|h\|^2,
$$
whence, due to the fact that $\nabla^2f(x)$ is symmetric positive semidefinite,
$$u^T\nabla^2f(x)v\leq [u^T\nabla^2f(x)u]^{1/2}[v^T\nabla^2f(x)v]^{1/2}\leq \alpha\|u\|\|v\|,$$ so that $(*_o)$ holds. Vice versa, assuming that $f$ obeys $(*_o)$, for $x,x+h\in V_o$ we have
$$
h^T[\nabla f(x+h)-\nabla f(x)]=\int_{0}^1h^T\nabla^2f(x+th)hdt\leq \int_{0}^1\alpha\|h\|^2dt=\alpha \|h\|^2.
$$
as required in $(i_o)$.
\par
$(2_o)\Leftrightarrow (*_o)$: Assuming that $(2_o)$ holds, for $x\in V_o$ and $u,v\in \bR^n$ we have
\[u^T\nabla^2f(x)v=\lim_{t\to +0} t^{-1}u^T[\nabla f(x+tv)-\nabla f(x)]\leq \varlimsup_{t\to+0}
t^{-1}\|u\|\|\nabla f(x+tv)-\nabla f(x)\|_*,
\] and the latter limit is $\leq \alpha\|u\|\|v\|$ by $(2_o)$. Vice versa, if $f$ obeys $(*_o)$, and $x\in V_o$, $x+h\in V_o$, we have
\begin{align*}
\|\nabla f(x+h)-\nabla f(x)\|_*&=\max_{u:\|u\|\leq1}u^T[\nabla f(x+h)-\nabla f(x)]=\max_{u:\|u\|\leq1}\int_{0}^1 u^T\nabla^2f(x+th)hdt\\
&\leq \max_{u:\|u\|\leq1}\int_{0}^1
\alpha \|u\|\|h\|dt=\alpha\|h\|,\\
\end{align*}
as required in $(2_0)$.
\par
\par
$(3_o)\Leftrightarrow (*_o)$: $(*_o)$ clearly implies $(3_o)$. On the other hand, when $x\in V_o$ and $h\in\bR^n$, for properly selected $C=C(x,h)<\infty$ and $\overline{t}=\aic{\overline{t,h}}{\overline{t}(x,h)}>0$ we have
$$
0\leq t\leq \overline{t}\Rightarrow
f(x+th)\geq f(x)+th^T\nabla f(x)+{t^2\over 2} h^T\nabla^2 f(x)h-Ct^3.
$$
It follows that when $f$ obeys $(3_o)$,
we have
$$f(x)+th^T\nabla f(x)+{t^2\over 2}h^T\nabla^2f(x)h-Ct^3\leq f(x)+th^T\nabla f(x)+{\alpha\over 2}\|th\|^2$$
for $0\leq t\leq \overline{t}$,  implying that $h^T\nabla^2f(x)h\leq\alpha\|h\|^2$, whence, $(*_o)$ holds true. \qed
%\section{Technical lemma}
\begin{lemma}\label{lemtech} Let $f(\cdot):\bR^n\to\bR_+$ be a norm such that the function $F(x)=f^2(x)$ is continuously differentiable and for certain $\alpha>0$ satisfies the relation
\[%begin{equation}\label{tleq1}
h^T[\nabla F(x+h)-\nabla F(x)]\leq \alpha F(h)\,\,\forall x,h\in\bR^n.
\]%end{equation}
Then \par
\item{(i)} For every $\epsilon\in(0,1)$ there exists a norm $f^\epsilon(\cdot)$ on $\bR^n$ such that
\[%begin{equation}\label{tleq2}
(1-\epsilon)f(\cdot)\leq f^\epsilon(\cdot)\leq (1+\epsilon)f(\cdot)
\]%end{equation}
and the function $F^\epsilon(x):=[f^\epsilon(x)]^2$ is C$^\infty$ outside of  the origin and satisfies the relations
\begin{align*}%\label{tleq3}
h^T[\nabla F^\epsilon(x+h)-\nabla F^\epsilon(x)]\leq (1+\epsilon)\alpha F^\epsilon(h)\quad &\forall (x,h\in\bR^n),\\
0\leq D^2F^\epsilon(x)[h,h]\leq \alpha F^\epsilon(h)\quad &\forall (x\neq0,h).
\end{align*}
Moreover, as $\epsilon\to+0$, functions $F^\epsilon$ converge to $F$ uniformly on compact subsets of $\bR^n$ along with their gradients.
\item{(ii)}  If, in addition, $f$ is an absolute norm, $f^\epsilon$ in {(i)} can be made absolute norm as well.
\end{lemma}
{\bf Proof.} {\em (i)} Let $S=\{x\in\bR^n:\|x\|_2=1\}$ and let $SO_n$ be the group of  orthogonal $n\times n$ matrices with determinant 1, both sets equipped with the natural structures of C$^\infty$ manifolds. Let also $\mu(\cdot)$ be the invariant measure on $SO_n$. In the sequel, for $g\in SO_n$ and function $\phi(\cdot)$ on $\bR^n$,  $\phi_g(\cdot)$ is given by $\phi_g(x)=\phi(gx)$.
\\{\bf 1$^o$.} Observe that for every $\epsilon\in(0,1)$ there exists $\delta(\epsilon)>0$ such that
$$
\|g-I_n\|_{2,2}\leq\delta<\delta(\epsilon)\;\Rightarrow\; (1-\epsilon)f(x)\leq f_g(x)\leq (1+\epsilon)f(x)\;\forall x
$$
where, same as above, $\|\cdot\|_{2,2}$ is the spectral norm. Indeed, for properly selected $C<\infty$ and any linear mapping $x\mapsto Hx:\bR^n\to\bR^n$  it holds
\[\|H\|_f\leq C\|H\|_{2,2}
\]
where $\|\cdot\|_f$ is the norm induced by the norm $f(\cdot)$ on $\bR^n$.
Consequently, for $g\in SO_n$, by the triangle inequality, one has for all $x$:
\begin{align*}
f(gx)&=f([g-I_n]x+x)\leq f(x)+\|g-I_n\|_ff(x)\leq (1+C\|g-I_n\|_{2,2})f(x),\\
f(gx)&=f([g-I_n]x+x)\geq f(x)-\|g-I_n\|_ff(x)\geq (1-C\|g-I_n\|_{2,2})f(x),
\end{align*}
so that it suffices to put $\delta(\epsilon)=\epsilon/C$.
\\
{\bf 2$^o$.} Given $\delta>0$, let $\theta(\cdot)$ be a nonnegative C$^\infty$ function on $SO_n$ such that $\theta(g)=0$ for $\|g-I_n\|_{2,2}\geq\delta$ and $\int_{SO_n}\theta_\delta(g)\mu(dg)=1$, and let
$$
F^{(\delta)}(x)=\int_{SO_n}F(gx)\theta_\delta(g)\mu(dg).
$$
Observe that function $F^{(\delta)}$ is continuous, homogeneous of degree 2 and positive outside of the origin, so that function $f^{(\delta)}(x)=\sqrt{F^{(\delta)}(x)}$ is positive outside of the origin and satisfies $f^{(\delta)}(\lambda x)=|\lambda|f(x)$. Besides this, \[f^{(\delta)}(x)=\sqrt{\int_{SO_n} f_g^2(x)\theta_\delta(g)\mu(dg)},\] so that $f^{(\delta)}$ is convex, since the functional $\Phi[\phi]=\sqrt{\int_{SO_n}\phi^2(g)\theta_\delta(g)\mu(dg)}$ on the space of continuous real-valued functions on $SO_n$ is convex and satisfies $\Phi[\psi]\geq \Phi[\phi]$ whenever $\psi(\cdot)\geq\phi(\cdot)\geq0$. The bottom line is that $f^{(\delta)}$ is a norm on $\bR^n$. In addition, taking into account that $\theta_\delta(g)\geq0$, $\int_{SO_n}\theta_\delta(g)\mu(dg)=1$, and $\theta_\delta(g)=0$ when $\|g-I_n\|_{2,2}\geq\delta$,
 by 1$^o$ we conclude that
\begin{equation}\label{tleq10}
\forall \epsilon\in(0,1),\; 0<\delta<\delta(\epsilon): \;(1-\epsilon)f(\cdot)\leq f^{(\delta)}(\cdot)\leq (1+\epsilon)f(\cdot).
\end{equation}
\\{\bf 3$^o$.} Since $F$ is continuously differentiable, we have
$$
\forall (x,h): \nabla F^{(\delta)}(x)=\int_{SO_n}g^T\nabla F(gx)\theta_\delta(g)\mu(dg),
$$
whence
\begin{align*}
h^T[\nabla F^{(\delta)}(x+h)-\nabla F^{(\delta)}(x)]&=\int_{SO_n} [gh]^T[\nabla F(gx+gh)-\nabla F(gx)]\theta_\delta(g)\mu(dg)\\
&\leq \int_{SO_n} \alpha F(gh)\theta_\delta(g)\mu(dg),
\end{align*}
that is,
\[%begin{equation}\label{tleq12}
\forall x,h:\; h^T[\nabla F^{(\delta)}(x+h)-\nabla F^{(\delta)}(x)]\leq \alpha F^{(\delta)}(h).
\]%end{equation}
{\bf 4$^o$.} Our goal is to prove that $F^{(\delta)}$ is C$^\infty$ outside of  the origin.
\par
Let $e_1,...,e_n$ be the canonic basis orths of $\bR^n$. For $u\in S$, let \[S(u)=\{v\in S:v^Tu>0\}=\{v\in S:\|v-u\|_2<\sqrt{2}\}.\]
We claim that there exists mapping $\chi:S(e_1)\to SO_n$ such that  $\chi(e_1)=I_n$ and $\chi(x)e_1=x$ for all $x\in S(e_1)$. Indeed, when $x\in S(e_1)$,
vectors $x,e_2,e_3,...,e_n$ are linearly independent. When applying to this sequence the Gram-Schmidt orthogonalization process, we obtain an orthonormal system $e_1(x)\equiv x,e_2(x),...,e_n(x)$ with vector-functions $e_i(x)$ which are C$^\infty$ on $S(e_1)$, and we can set $\chi(x)=[e_1(x),e_2(x),...,e_n(x)]$.\par
We are now ready to show that $F^{(\delta)}$ is C$^\infty$ outside of  the origin when $\epsilon\in(0,1)$ and $0<\delta<\delta(\epsilon)$. Since $F^{(\delta)}$ is homogeneous of degree 2 and positive outside of   the origin, it suffices to verify that the restriction of $F^{(\delta)}$ onto the sphere $S$ is C$^\infty$. To justify the latter claim, let $\delta>0$, $\bar{x}\in S$, and let $\bar{g}\in SO_n$ be such that $\bar{g}e_1=\bar{x}$. Then for all $x\in S(\bar{x})$
\[ \bar{g}^{-1}x\in S(e_1)\Rightarrow \bar{g}^{-1}x=\chi(\bar{g}^{-1}x)e_1\Rightarrow
x=\bar{g}\chi(\bar{g}^{-1}x)e_1,\]
whence for such $x$ and $h=g\bar{g}\chi(\bar{g}^{-1}x)$ we get
\[
%\begin{array}{l}
%x\in S(\bar{x})\Rightarrow
F^{(\delta)}(x)=\int_{SO_n}F(g\bar{g}\chi(\bar{g}^{-1}x)e_1)\theta_\delta(g)\mu(dg)\\
=\int_{SO_n}F(he_1)\underbrace{\theta_\delta(h\chi^{-1}(\bar{g}^{-1}x)\bar{g}^{-1})}_{\theta_{\bar{x}}(x,h))}\mu(dh)
\]
due to invariance of $\mu(\cdot)$.
Because $\theta(\cdot)$ is C$^\infty$ on $SO_n$ and $\chi(z)$ is C$^\infty$ in $S(e_1)$, function $\theta_{\bar{x}}(\cdot,h)$ is C$^\infty$ in the first argument in the neighbourhood $S(\bar{x})$ with derivatives which are continuous in $x$ and $h$, implying that $F^{(\delta)}(x)$ is C$^\infty$ on $S$.
\\
{\bf 5$^o$.} The bottom line of the above is that for every $\epsilon\in(0,1)$, for properly selected $\delta(\epsilon)>0$ and all $\delta\in(0,\delta(\epsilon))$ function  $F^{(\delta)}(\cdot)$ is the square of some norm $f^{(\delta)}(\cdot)$ satisfying
$$
(1-\epsilon)f(\cdot)\leq f^{(\delta)}(\cdot)\leq (1+\epsilon)f(\cdot),
$$
$F^{(\delta)}$ is C$^\infty$ outside of the origin, and
\begin{equation}\label{out1}
\forall x,h: h^T[\nabla F^{(\delta)}(x+h)-\nabla F^{(\delta)}(x)]\leq \alpha F^{(\delta)}(h).
\end{equation}
In particular, when $x\neq 0$, we have
$$
0\leq h^T[\nabla F^{(\delta)}(x+th)-\nabla F^{(\delta)}(x)]\leq t^{-1}\alpha F^{(\delta)}(th)=t\alpha F^{(\delta)}(h),
$$
implying that
\begin{equation}
  \label{out2}
D^2F^{(\delta)}(x)[h,h]=\lim_{t\to+0}t^{-1}h^T[\nabla F^{(\delta)}(x+th)-\nabla F^{(\delta)}(x)]\leq \alpha F^{(\delta)}(h)\,\forall (x\neq0,h).\end{equation}
Thus, given $\epsilon\in(0,1)$, the function $F^\epsilon:=F^{(\delta(\epsilon)/2)}$ satisfies all requirements of {\em (i)}. {\em (i)} is proved.
\\
{\bf 6$^o$.} It remains to prove {\em (ii)}. To this end, assume that $f$ is an absolute norm. Let ${\cal E}_n$ be the multiplicative group comprised by the $2^n$ diagonal matrices with diagonal entries $\pm1$. Given $\epsilon\in(0,1)$ and assuming that $\delta\in(0,\delta(\epsilon))$, let us set
$$
F^{\delta,E}(x):=F^{(\delta)}(Ex)=[f^{(\delta)}(Ex)]^2,\,E\in {\cal E}_n.
$$
Functions $F^{\delta,E}(x)$ are C$^\infty$ outside of the origin functions, while due to (\ref{out1}) and (\ref{out2}) we have
\begin{align*}
\forall x,h:\;\; h^T[\nabla F^{\delta,E}(x+h)-\nabla F^{\delta,E}(x)]&\leq \alpha F^{\delta,E}(h),\\
\forall x\neq0,h:\;\;D^2F^{\delta,E}(x)[h,h]&\leq \alpha F^{\delta,E}(h).
\end{align*}
We conclude that the function
$$
\overline{F}^{\delta}(x)=2^{-n}\sum_{E\in {\cal E}_n} F^{\delta, E}(x)
$$
is C$^\infty$ outside of  the origin and satisfies
\begin{align*}
\forall x,h: \;\;h^T[\nabla \overline{F}^{\delta}(x+h)-\nabla \overline{F}^{\delta}(x)]&\leq \alpha \overline{F}^{\delta}(h),\\
\forall x\neq0,h:\;\;D^2\overline{F}^\delta(x)[h,h]&\leq \alpha \overline{F}^\delta(h).\end{align*}
Besides this,
\[
\left(\overline{F}^\delta(x)\right)^{1/2}=\left(2^{-n}\sum_{E\in{\cal E}_n}[f^{(\delta)}(Ex)]^2\right)^{1/2},
\]
so that $\left(\overline{F}^\delta(x)\right)^{1/2}$ is a norm along with the functions $f^{(\delta)}(Ex)$. Furthermore, by construction, $\overline{F}^\delta(Ex)=\overline{F}^\delta(x)$ for all $x$, so that $\left(\overline{F}^\delta(x)\right)^{1/2}$ is an absolute norm. Finally, by (\ref{tleq10}) and due to $0<\delta<\delta(\epsilon)$ we have
\[
(1-\epsilon)^2F(\cdot)\leq [f^{\delta}(\cdot)]^2\leq (1+\epsilon)^2F(\cdot),
\]
and because $f(\cdot)$ is an absolute norm, it follows that
\[
(1-\epsilon)^2F(x)\leq [f^{{\delta}}(Ex)]^2\leq (1+\epsilon)^2F(x)\;\;\forall x,\forall E\in{\cal E}_n.
\]
Hence,
\[
(1-\epsilon)^2F(\cdot)\leq \overline{F}^{\delta}(\cdot)\leq (1+\epsilon)^2F(\cdot),
\]
and thus
$$
(1-\epsilon)f(\cdot)\leq \left({\overline{F}^{\delta}(\cdot)}\right)^{1/2}\leq (1+\epsilon)f(\cdot).
$$
The bottom line is that to meet all requirements of {\em (ii)}, it suffices to set $F^{\epsilon}=\overline{F}^{\delta(\epsilon)/2}$.
\\
{\bf 7$^o$.} Finally, we can ensure that $\delta(\epsilon)\to+0$ as $\epsilon\to+0$, and in this case, as is immediately seen from the above, functions $F^\epsilon$ converge, along with their gradients, to $F$ uniformly on bounded subsets of $\bR^n$ as $\epsilon\to +0$. \qed

\end{document}